\documentclass{article}
\usepackage{heron2e}
\usepackage{times} 
\title{Weak C* Hopf symmetry\footnote{Talk at the XXI International 
Colloquium on Group Theoretical Methods in Physics, Goslar 1996}}
\author{Karl-Henning Rehren}
\address{II. Institut f\"ur Theoretische Physik, Universit\"at Hamburg, \\ 
Luruper Chaussee 149, D-22761 Hamburg, Germany, \\
email: rehren@x4u2.desy.de}

\def\sitem#1{{\small (#1)}} \def\sbitem#1{\par{\small\bf (#1)}}
\def\aa{{\cal A}} \def\bb{{\cal B}} \def\oo{{\cal O}}
\def\sig{{\sigma}}  \def\eps{\varepsilon}  \def\rho{\varrho}  \def\id{{id}}
\def\plus{{_\oplus}}   \def\red{_{\rm red}}  \def\reg{_{\rm reg}}
\def\comp{{\scriptstyle \circ}}  \def\inv{^{-1}}     
\def\arr#1{\colon#1\!\rightarrow\!} \def\act{\triangleright}
\def\VVert{\hbox{\hbox{$\vert$}\kern-1pt\hbox{$\vert$}\kern-1pt\hbox{$\vert$}}}
\def\solv{\vbox{\vskip-7pt
\hbox{\lower8pt\hbox{\kern.9pt$\scriptscriptstyle\bullet$}}\hbox{$v$}}{}}
\def\solw{\vbox{\vskip-7pt
\hbox{\lower8pt\hbox{\kern2.2pt$\scriptscriptstyle\bullet$}}\hbox{$w$}}{}}
\def\opw{\vbox{\vskip-7pt
\hbox{\lower8pt\hbox{\kern2.2pt$\scriptscriptstyle\circ$}}\hbox{$w$}}{}}
\def\opv{\vbox{\vskip-7pt
\hbox{\lower8pt\hbox{\kern.9pt$\scriptscriptstyle\circ$}}\hbox{$v$}}{}}
\def\iotasubset{\mathrel{\vbox{\vskip-2pt
      \hbox{\kern3.5pt \lower5pt\hbox{$\scriptstyle\iota$}}
      \hbox{\hbox{$\subset$}\kern-5pt\lower2.4pt\hbox{$\rightarrow$}}}}}
\def\hatQ{\vbox{\vskip-2pt\hbox{$\widehat Q$}}}
\def\semi{\setbox0=\hbox{$\times$}
     \hbox to0pt{\kern0.78\wd0 \vrule width.05ex height0.8\ht0\hss}\box0}
\def\eins{{\mathchoice {\rm 1\mskip-4mu l} {\rm 1\mskip-4mu l}
  {\rm 1\mskip-4.5mu l} {\rm 1\mskip-5mu l}}}
\def\CC{{\mathchoice {\setbox0=\hbox{$\displaystyle\rm C$}\hbox{\hbox
  to0pt{\kern0.4\wd0\vrule height0.9\ht0\hss}\box0}}
  {\setbox0=\hbox{$\textstyle\rm C$}\hbox{\hbox
  to0pt{\kern0.4\wd0\vrule height0.9\ht0\hss}\box0}}
  {\setbox0=\hbox{$\scriptstyle\rm C$}\hbox{\hbox
  to0pt{\kern0.4\wd0\vrule height0.9\ht0\hss}\box0}}
  {\setbox0=\hbox{$\scriptscriptstyle\rm C$}\hbox{\hbox
  to0pt{\kern0.4\wd0\vrule height0.9\ht0\hss}\box0}}}}
\begin{document}
\maketitle
\begin{abstract} 
Weak C* Hopf algebras can act as global symmetries in low-dimensional
quantum field theories, when braid group statistics prevents ordinary group
symmetries. Charged fields transform linearly in finite multiplets, and the 
observables are precisely the gauge invariants. Possibilities to construct 
field algebras with weak C* Hopf symmetry from a given theory of local 
observables are discussed.
\end{abstract}

\section{Introduction}
On several occasions at the Arnold Sommerfeld Institute at Clausthal 
(and elsewhere) I have addressed the field problem and the symmetry problem 
in the theory of superselection sectors of low-dimensional models of quantized 
fields \cite{CLZi,CLZii}. 

The field problem raises the question whether the algebra of physical 
observables can be embedded into a field algebra containing sufficiently many 
charged fields to generate from the vacuum all superselection sectors of the 
observables. 
The identification of a sensible class of gauge symmetries in a given model 
with braid group statistics (which precludes symmetry groups) constitutes 
the symmetry problem.

The field problem can always be solved by a canonical construction of 
intertwining field operators, which was introduced in \cite{FRS} and called 
the ``reduced field bundle'' (with commutation relations in the form of an 
exchange algebra). 
Yet, I have myself rejected this answer as a useful candidate for the
ambient field algebra. 
The reason was the apparent absence of a decent symmetry under which the 
reduced field bundle transforms, with the observables as the invariants 
(the ``gauge principle''). 

A main point in this report is the observation that the reduced field bundle
does possess a symmetry which I have previously overlooked. 
In fact, it is invariant under the action of a weak C* Hopf algebra. 
Weak C* Hopf algebras \cite{BSz,NW} (and related objects \cite{Y}) are 
comparatively conservative generalizations of finite symmetry groups. 
The coproduct is non-commutative and does not map the unit operator onto the 
unit operator. 
Unlike quantum groups \cite{Dr}, however, it preserves the $*$-structure, 
and unlike quasi quantum groups \cite{Sch}, it is coassociative. 
Unlike the combinatorial concept of a paragroup \cite{O}, and unlike my 
personal attempts with C* symmetries \cite{CLZii}, its action on the field 
algebra is determined by a linear transformation law on finite tensor 
multiplets of charged fields.

Weak C* Hopf symmetries are by no means restrained to the reduced 
field bundle. 
They are rather the natural symmetry concept associated with
finite index subfactors of depth 2. 
If such a subfactor is irreducible, then the symmetry is in fact a C* Hopf 
algebra; so its ``weakness'' is precisely due to reducibility. 
The depth 2 condition has rather a technical meaning (``all irreducible 
components of the subfactor $A \subset B_2$ are already contained in 
$A \subset B$''), and in the case of irreducible subfactors becomes a simple 
structural property (``$A' \cap B_2$ is a factor''). 
Here $B_2$ is the second Jones extension of $A \subset B$. 

Weak C* Hopf algebras possess unitary matrix representations.
The size of the corresponding multiplets (of charged fields) is, in 
contrast to the true Hopf case, different from the quantum dimension of
the corresponding bimodule (statistical dimension of the associated
superselection sector). 
Consequently, the obstruction against non-integer quantum (or statistical) 
dimensions which seemed always inherent to C* symmetries, is absent for 
weak C* Hopf symmetries.

\section{Weak C* Hopf algebras}
The following defining system of axioms for weak C* Hopf algebras is 
due to ref.\ \cite{BSz}, to which we also refer for further details. 
Our emphasis here is on the departure from true C* Hopf algebras. 

\sbitem A A finite dimensional weak C* Hopf algebra is a C* 
algebra with $\eins$ (hence a direct sum of matrix rings) with the additional 
structures coproduct, counit and antipode.

The coproduct is a coassociative $*$-homomorphism $\Delta \arr Q Q \otimes Q$. 
The counit is a positive linear map $\eps \arr Q \CC$ and satisfies the 
compatibility condition with the coproduct: $(\eps\otimes\id) \comp \Delta = 
\id = (\id\otimes\eps) \comp \Delta$. 
The antipode is a complex-linear anti-homomorphism and anti-cohomomorphism 
$S \arr QQ$ (i.e., it reverts the order of the product and of the coproduct),
and is inverted by the $*$-structure: $S\inv(q) = S(q^*)^*$.

\sbitem B The three axioms in the left column hold, as opposed to the 
corresponding three stronger axioms for true C* Hopf algebras in the 
right column.
$$ \hbox{Weak:}\left\{ 
\matrix{ \Delta(\eins) \equiv \eins_{(1)} \otimes \eins_{(2)} = 
\hbox{Projection} \cr \eps(qp) = \eps(q\eins_{(1)})\cdot\eps(\eins_{(2)}p) \cr
S(q_{(1)})q_{(2)} \otimes q_{(3)} = (\eins \otimes q) \cdot \Delta(\eins) \cr}
\right. \qquad \hbox{True:}\left\{ \matrix{ \Delta(\eins) = \eins 
\otimes \eins \cr \eps(qp) = \eps(q)\cdot\eps(p) \cr S(q_{(1)})q_{(2)} = 
\eps(q)\eins \cr} \right. $$
(We use the shorthand notation $q_{(1)} \otimes q_{(2)}$ for the expansion of 
$\Delta(q)$.) 
Along with the three weak axioms, either of the three stronger axioms implies 
the other two. Therefore these three axioms for true C* Hopf algebras 
cannot be independently relaxed.

\sbitem C The dual $\hatQ$ is defined by the linear maps $\{\hat q\arr Q 
\CC\}$. 
The structure data of $Q$ are canonically dualized by the pairing, and by 
$\langle \hat q{}^*,q \rangle = \overline{\langle \hat q,S(q)^* \rangle}$ 
the dual is given a $*$-structure. 
$\hatQ$ with data $(\hat\eins = \eps,\hat\Delta,\hat\eps,\hat S)$ and $*$ is 
again a weak C* Hopf algebra.

\sbitem D A (left) action of a weak C* Hopf algebra $Q$ on a C* 
algebra $B$
is a unital algebra homomorphism from $Q$ into the linear maps of $B$ into $B$,
denoted by $b \mapsto (q \act b)$, satisfying 
$$ \matrix{ q \act (b\cdot c) = (q_{(1)} \act b)\cdot(q_{(2)} \act c),
\cr \left( q^* \act b^* \right)^* = S\inv(q) \act b \cr} \qquad 
\left(\matrix{ q \in Q \cr b,c \in B \cr}\right) 
\hbox{\raise-3mm\hbox{.}} $$
An element $a$ of $B$ is called invariant under the action, if
$q \act a = q_{(1)}S(q_{(2)}) \act a$.
The invariant elements form a $*$-subalgebra $A \equiv B^Q \subset B$.

\section{Subfactors of depth 2}
We consider a pair $A$, $B$ of von Neumann factors of type III, along 
with an injective unit-preserv\-ing homomorphism $\iota\arr AB$. $A$ may be 
thought of as a subfactor of $B$, and $\iota$ the inclusion map. 

\sbitem A \cite{Lo} The subfactor $A \iotasubset B$ has finite index if 
and only if there is a homomorphism $\bar\iota\arr BA$ (the conjugate) and a 
``standard'' pair of isometries $w \in A$, $v \in B$ satisfying 
$$ \matrix{ w a = \rho(a) w \qquad (a \in A), \cr
w^*\bar\iota(v) = \lambda^{-{1 \over 2}} \cdot \eins_A, \cr}
\qquad\qquad \matrix{ v b = \gamma(b) v \qquad (b \in B), \cr
v^*\iota(w) = \lambda^{-{1 \over 2}} \cdot \eins_B, \cr} $$
where $\gamma = \iota\comp\bar\iota$ and $\rho = \bar\iota\comp\iota$ are the 
canonical and dual canonical endomorphisms of $B$ and $A$, respectively, 
associated with the subfactor. $\lambda > 1$ is called the index of the 
subfactor associated with the pair $w,v$. In the following we assume $w$ and 
$v$ to minimize the index. 

\sbitem B \cite{LH} A ``canonical triple'' $(\rho,w,x)$ consists of an 
endomorphism $\rho$ of $A$ and a pair of isometries $w,x$ in $A$ satisfying 
the relations ($\lambda^2$ being the index of $\rho(A) \subset A$)
$$ \matrix{ w a = \rho(a) w, \qquad x \rho(a) = \rho^2(a) x \qquad 
(a \in A); \cr x x = \rho(x) x, \quad x x^* = \rho(x^*) x, \quad
w^* x = \lambda^{-{1 \over 2}}\eins_A = \rho(w^*) x. \cr} $$
These relations make certain that $\rho$ is a canonical endomorphism, i.e.,
one can draw a ``conjugate square root'' $\rho = \bar\iota\comp\iota$ 
with $w,v$ as in \sitem A and $x = \bar\iota(v)$.

Every such triple characterizes, up to isomorphism, an ambient algebra $B$
such that 
$$ \rho(A) \subset A_1 \subset A \iotasubset B $$ 
is a sequence of Jones extensions of index $\lambda$. 
Here $A_1$ is the intermediate algebra generated by $\rho(A)$ and $x$,
from which $B$ is obtained by the Jones construction. 
(An additional non-redundancy condition on the triple prevents $B$ from having 
a center. We shall tacitly always assume that $B$ is a factor.)
The minimal conditional expectation from $B$ onto $A$ is $\mu = 
w^*\bar\iota(\,\cdot\,)w$.

\sbitem C \cite{NW} A dual pair of weak C* Hopf algebras is associated 
with every reducible subfactor of finite index $\lambda$ and of depth 2. 
As algebras, these are the relative commutants $Q := \bar\iota\comp\iota(A)' 
\cap A$ and $\hatQ := \iota\comp\bar\iota(B)' \cap B$. 
They are direct sums of matrix rings corresponding to the subsectors of 
$\bar\iota\comp\iota$ resp.\ $\iota\comp\bar\iota$ and their multiplicities, 
and have the same finite dimension. 

Let $z$ be the positive central element in the relative commutant 
$\iota(A)' \cap B$ with value $\sqrt{n_s/d_s}$ on each of the minimal central 
projections $p_s$ of the relative commutant, where $d_s^2$ is the index of 
the corresponding reduced subfactor, and $n_s$ is its multiplicity. 
Let $\bar z$ in the relative commutant $\bar\iota(B)' \cap A$ be defined 
analogously, and put $\solv := zv \equiv \iota(\bar z)v$,  
$\solw := \bar z w \equiv \bar\iota(z)w$, and
$\opv := z\inv v \equiv \iota(\bar z\inv)v$,  
$\opw := \bar z\inv w \equiv \bar\iota(z\inv)w$.

The algebras $Q$ and $\hatQ$ are put into duality by the nondegenerate pairing
$$ \langle \hat q,q \rangle := \lambda \cdot w^*\bar\iota(\solv)^*
q \bar\iota(\hat q) \solw w \equiv \lambda \cdot v^*\iota(\solw)^*
\hat q \iota(q) \solv v. $$
It induces a coproduct in $Q$ (as a linear map into $Q \otimes Q$) from the 
product in $\hatQ$.
This coproduct is in fact a $*$-homomorphism if and only if the depth is 2,
i.e., if and only if every subsector of $\iota\comp\bar\iota\comp\iota$ is 
already contained in $\iota$.

The counit on $Q$ is the pairing with the unit in $\hatQ$. 
The antipode on $Q$ is given by the two equivalent definitions
$$ S(q) = \lambda \cdot \bar\iota\left(\iota\left[\opw^*\bar\iota(\solv)^*q 
\right]\solv\right)\opw \equiv \lambda \cdot \solw^* 
\bar\iota\left(\opv^*\iota\left[q\bar\iota(\opv)\solw \right]\right) $$
and is involutive up to a non-unitary conjugation:
$$ S(S(q)) = [\bar z\inv\bar\iota(z)]^2 \cdot q \cdot 
[\bar z\inv\bar\iota(z)]^{-2}. $$

The dual structures on $\hatQ$ are given by the replacements $(A,w,\iota) 
\leftrightarrow (B,v,\bar\iota)$. 
All statements here and in the remainder of this chapter are perfectly 
symmetric under this duality.
Our emphasis will, however, be on those statements which pertain to the natural
interpretation $A=$ fixpoints of $B$ under a symmetry.
\vskip 1mm

{\bf Proposition:} {\sl With the above definitions, 
$Q$ and $\hatQ$ are a dual pair of finite dimensional weak C* Hopf algebras.
The dual formulae
$$ q \act b := \lambda^{1 \over 2} \cdot \opv^*\iota\left(q\bar\iota(b)
\solw\right) \qquad \hbox{and} \qquad  \hat q \act a := \lambda^{1 \over 2} 
\cdot \opw^*\bar\iota\left(\hat q\iota(a)\solv\right) $$
define an action of $Q$ on the algebra $B$ with invariants $B^Q = \iota(A)$,
and an action of $\hatQ$ on $A$ such that $B$ is isomorphic with the
crossed product of $A$ by $\hatQ$.} 
\vskip 1mm

The last statement of the proposition means that $\hatQ$ and $A$ are 
embedded as subalgebras into $B$ with relations 
$$ \hat q \cdot \iota(a) = \iota(\hat q_{(1)}\act a) \cdot \hat q_{(2)} \qquad
\hbox{or equivalently} \qquad \iota(\hat q \act a) = \hat q_{(1)} \cdot 
\iota(a) \cdot \hat S(\hat q_{(2)}), $$ 
and together generate $B$. 

\sbitem D $B$ is spanned by elements $\Gamma_e\iota(a)$ where $a \in A$ and 
$\Gamma_e \in B$ are intertwiners $\iota\comp\rho_c \to \iota$ associated with 
every irreducible subsector (``charge'') $\rho_c \prec \bar\iota\comp\iota$.
Under the action of $Q$ the operators $\Gamma_e$ of fixed charge transform as 
finite multiplets by matrix multiplication of the corresponding matrix ring 
of $Q$. The dimension of each multiplet equals the multiplicity of $\rho_c$
in $\bar\iota\comp\iota$.

Acting on $B$ with the element $\opw\solw^*$ of $Q$ (which is related to the 
Haar measure) annihilates all multiplets of nontrivial charge and averages 
over those of trivial charge. The effect is precisely the minimal conditional 
expectation $\mu$ onto $A$.

\sbitem E The action of $Q$ on $B$ is partly inner in the following sense. 
The relative commutant $\iota(A)' \cap B$ is mapped by $\bar\iota$ onto a 
subalgebra $P$ of $Q$ such that $S(P) \equiv \bar P = \bar\iota(B)' \cap A$. 
$P$ and $\bar P$ are two commuting subalgebras of $Q$, and $P \!\cdot\! \bar P$
is a weak C* Hopf subalgebra of $Q$. One has $\Delta(p) = \Delta(\eins)
(p \otimes \eins)\Delta(\eins)$ for $p \in P$ \cite{BSz}
and $\eins_{(1)}\bar\iota(b)S(\eins_{(2)}) = \bar\iota(\eins \act b) = 
\bar\iota(b)$, implying 
$$ \bar\iota(p \bar p \act b) 
= p \cdot \bar\iota(b) \cdot S(\bar p) \qquad (p \in P, \bar p \in \bar P). $$ 
Since $P = S(\bar P) \subset \bar\iota(B)$ it follows that the action of 
$P \!\cdot\! \bar P$ on $B$ is unitarily implemented by elements of $B$. 
This feature is not surprising since also for group actions it is well known 
that an inner action produces a reducible inclusion. 

\sbitem F \cite{NW} The action of $Q$ on $B$ preserves the subalgebra 
$\hatQ$, and $\hatQ$ preserves the subalgebra $Q$ of $A$. 
The action of $Q$ on $\hatQ$ satisfies (and is determined by) the rule
$\langle q \act \hat q,p \rangle = \langle \hat q,pq \rangle$, and vice versa.

\section{Application to the field problem}
\indent
\sbitem A A quantum field theory assigns to every bounded region $\oo$ in 
space-time
a weakly closed operator algebra $\aa(\oo)$ generated by the fields (in a 
vacuum representation) localized in $\oo$. 
The resulting isotonous net of algebras determines the theory, even without 
knowledge of the underlying fields.
Under standard assumptions of covariance and spectrum condition, the local
algebras are hyperfinite type III factors.

The theory is local if the local algebras associated with two regions 
at space-like distance commute with each other.
A theory of observables has to be local.

The following results can be found in full detail in \cite{LR}.

\sbitem B A field extension of a theory of observables $\aa$ is a 
(relatively local, but possibly nonlocal) isotonous net $\bb$ such that for 
every region, $\aa(\oo)$ is a subfactor of $\bb(\oo)$. 
Assuming the existence of a conditional expectation $\mu\arr\bb\aa$ which 
preserves the localization and leaves the vacuum state invariant (an unbroken 
global symmetry in the broadest sense), the dual canonical endomorphism $\rho$ 
associated with a single local subfactor $\aa(\oo) \subset \bb(\oo)$ extends 
to a covariant endomorphism $\rho$ of the entire net of observables $\aa$. 
The latter is trivial on observables at space-like distance from $\oo$,
and is hence a DHR endomorphism \cite{DHR} localized in $\oo$. 
$\rho$ is equivalent to the representation of $\aa$ on the vacuum Hilbert 
space of $\bb$, and its subsectors are precisely those superselection 
charges of $\aa$ for which there are charged fields within the ambient net 
$\bb$.

\sbitem C Conversely, every ``DHR canonical triple'' $(\rho,w,x)$ 
determines, in terms of observable data, a field extension $\bb$ with an 
unbroken global symmetry. 
Here $\rho$ is a DHR endomorphism of $\aa$ localized in some region $\oo$, 
$w$ and $x$ are isometries in $\aa(\oo)$, and the algebraic relations as in 
3\sitem B hold with $a \in \aa$. 

The algebra $B \equiv \bb(\oo)$ is constructed from the triple as in 3\sitem B 
with $A \equiv \aa(\oo)$, and the other local algebras are obtained from it by 
Poincar\'e or conformal covariance. 
It is a nontrivial result about this construction that the structural 
properties of the ``germinal'' local subfactor $\aa(\oo) \subset \bb(\oo)$ 
propagate to every other local subfactor. 
In the sequel, statements about the net are understood to hold for every local 
subfactor. 

The simple eigenvalue condition $\eps_\rho x = x$ on the statistics operator 
$\eps_\rho$ of $\rho$ decides whether the net $\bb$ resulting from the triple 
is local or not.
 
\sbitem D A ``cheap'' way to obtain DHR canonical triples is the following. 
For any DHR endomorphism $\sig$ of $\aa$ (with finite statistics) there is a 
standard pair (cf.\ 3\sitem A) of isometries $w$ and $\bar w$ in $\aa$ such 
that 
$wa = \sig\bar\sig(a)w$ and $\bar wa = \bar\sig\sig(a)\bar w$ for $a \in \aa$. 
Thus $(\rho=\sig\comp\bar\sig,w,x=\sig(\bar w))$ is a DHR canonical triple.
Therefore, every DHR sector $\sig$ defines a corresponding field 
extension $\bb_\sig$.
If $\sig$ is localized in $\oo$, then $B \equiv \bb_\sig(\oo)$ is the Jones 
extension of $A \equiv \aa(\oo)$ by its subfactor $A_1 \equiv \sig(\aa(\oo))$.
The local subfactors have depth 2 if and only if all subsectors of 
$\sig\bar\sig\sig$ are already contained in $\sig$. 

By the eigenvalue condition (cf.\ \sitem C), the nets $\bb_\sig$ are always 
nonlocal (unless $\sig$ is an automorphism and consequently $\bb_\sig = \aa$). 
Since this statement is not in \cite{LR}, we provide the argument here.
One applies the standard left-inverse $\phi(a) = \bar w^*\bar\sig(a)\bar w$
of $\sig$ to the eigenvalue condition. But $\phi(x)=\bar w$ is an isometry 
while, by the statistics calculus \cite{FRS,DHR}, 
$\phi(\eps_{\sig\bar\sig}\sig(\bar w))$ differs from an isometry by the 
statistics parameter of $\sig$ which is $1/d(\sig)$ times a unitary. 
Therefore, equality can hold only if $d(\sig)=1$.

\sbitem E This ``cheap''
construction does not provide all DHR canonical triples.
Notably an irreducible extension with an outer action of a compact gauge group 
such that the observables are the fixed points cannot be of this type.
E.g., provided all DHR sectors of $\aa$ have permutation group statistics, the 
Doplicher-Roberts construction \cite{DR} determines a graded local extension 
with a compact gauge group. 
Its dual canonical endomorphism $\rho = \rho\reg$ contains every irreducible 
DHR sector with multiplicity $n_s$ equal to its statistical dimension $d_s$. 
The isometry $x$ encodes the Clebsch-Gordan coefficients of the gauge group.

\section{Putting things together}
In a quantum field theory with superselection sectors, it is desirable to
have a field algebra of charged fields which generate all charged sectors
from the vacuum, and a gauge symmetry acting on the fields with the 
observables as fixpoints. 
In order to be of practical use, the transformation law for the fields 
should be sufficiently simple and concrete.
These demands favour field extensions with weak C* Hopf symmetry which have 
transformation laws with finite multiplets, while in low dimensions it is 
in general inconsistent to ask for true C* Hopf symmetry.

\sbitem A In order to find field extensions of a given (rational) local 
quantum field theory of observables $\aa$ which have a weak C* Hopf 
symmetry, one has to look for DHR canonical triples $(\rho,w,x)$ such that 
the resulting local subfactors have depth 2.

\sbitem B An immediate possibility in rational models is to choose $\sig = 
\sig\plus$ the direct sum of all irreducible DHR sectors of 
$\aa$ with multiplicity one, and to proceed as in 4\sitem D. 
By construction and since $\sig\plus$ is self-conjugate, the resulting 
local subfactor $\aa(\oo) \subset \bb(\oo)$ is isomorphic to 
$\sig\plus(\aa(\oo)) \subset \aa(\oo)$. 
An explicit unitary equivalence shows \cite{N} that the same holds 
for the reduced field bundle extension $\bb\red$.
Therefore, the reduced field bundle equals the extension of $\aa$ by 
$\sig\plus$. 
It has depth 2 and a self-dual weak C* Hopf symmetry $Q$.
As algebras, $Q$ and $\hatQ$ are both isomorphic to the relative commutant 
$\sig\plus^2(\aa)' \cap \aa$. 
The latter is a sum of matrix rings $M_s \cong \hbox{Mat}_{N_s}(\CC)$ 
labelled by the DHR sectors (``charges'') of $\aa$, with $N_s$ equal to their 
multiplicities within $\sig\plus^2$.

The irreducible multiplets $\Gamma_e\iota(a)$ (cf.\ 3\sitem D) coincide with 
the charged operators $F_e(a)$ spanning the reduced field bundle as defined 
in \cite{FRS}. 
Averaging over the symmetry by the conditional expectation $\mu$ 
(cf.\ 3\sitem{B,D}) yields the invariant elements $\sum_e F_e(a)$ where the 
sum extends over all edges of trivial charge. 
These are precisely the observables $\iota(a)$ on the extended Hilbert space.

\sbitem C If all DHR sectors of $\aa$ have integer statistical dimension, 
another natural choice is $\sig = \sig\reg \equiv \rho\reg$ ($n_s = d_s$, 
$z = \eins = \bar z$).
This choice is possible whenever the observables are given as the 
fixed points of another net $\bb$ under a finite gauge group, e.g., the 
graded local Doplicher-Roberts extension in the case of permutation group 
statistics, cf.\ 4\sitem E. 

In the diagram below, the first row is a Jones sequence due to the
gauge symmetry, and so is the second row by definition of $\bb\reg$. 
Since an alternating subsequence of a Jones sequence is again a Jones 
sequence, the two obvious vertical equalities imply the third:
$$ \matrix{A_2 = \rho\reg(A) & \subset & A_1 & \subset & A & \subset & B 
& \subset & B_1 = B \semi G \cr 
\qquad\quad\VVert &&&& \VVert &&&& \Vert\qquad\quad\;\; \cr
\qquad\;\;\sig\reg(A) && \subset && A && \iotasubset && B\reg\qquad\quad \cr} 
\qquad\quad 
\pmatrix{A = \aa(\oo) \cr B = \bb(\oo) \cr} \hbox{\raise-3mm\hbox{.}} $$
It follows that the extension $\bb\reg$ of $\aa$ by $\sig\reg(\aa)$ equals 
the extension of $\bb$ by $\aa$, which in turn equals the crossed product of 
$\bb$ by its gauge group $G$.

\sbitem D If all sectors of $\aa$ are simple, then $\sig\reg$ coincides 
with $\sig\plus$. 
There is an anyonic field extension $\bb$ with an abelian symmetry group such 
that the dual canonical endomorphism is $\rho\reg = \sig\reg = \sig\plus$. 
By combination of \sitem B and \sitem C, the reduced field bundle is a crossed 
product of the anyonic extension by its abelian gauge group.

\sbitem E In the general case with integer statistical dimensions, the Jones 
extension $\bb\reg$ of $\aa$ by $\sig\reg(\aa)$ is to the same extent
``larger'' than $\bb\red$ as $\sig\reg(\aa)$ is ``smaller'' than 
$\sig\plus(\aa)$. 
Clearly, one has no inclusion as algebras but rather a compression by an 
appropriate projection which selects $\sig\plus \prec \sig\reg$.

\sbitem F By combination of \sitem C and \sitem E, if the sectors of $\aa$ 
have 
permutation group statistics, then the reduced field bundle is a compression 
of the crossed product of the Doplicher-Roberts graded local field algebra by 
its gauge group.
 
\section{Discussion}
The paragroups assocciated with subfactors of finite index and depth 2 are 
weak C* Hopf algebras. 
Weak C* Hopf algebras therefore reconcile the generalized symmetry notions 
due to Ocneanu and due to Kac and Drinfel'd. 
They arise as symmetries in quantum field theoretical models whenever the 
local subfactors ``observables $\subset$ fields'' have depth 2 (the condition 
of finite index may presumably be relaxed).

Field extensions with global weak C* Hopf symmetries are encoded 
in terms of the observables by depth 2 DHR canonical triples $(\rho,w,x)$.
An obvious and not very subtle class of such extensions, including the reduced 
field bundle, is described in 4\sitem D. 
Relations between various such extensions are clarified. 
Notably if the observables are the fixed points under a finite nonabelian 
gauge group acting on a given field net, it is made clear in which sense the 
reduced field bundle exceeds the given field net (it corresponds to its 
crossed product by the gauge group), and in which sense it is smaller than 
the former (it is a compression which removes the gauge multiplicities).

For a given theory there may exist other (and more ``economic'' as far as the 
field problem is concerned) DHR canonical triples which are not of the 
form described in 4\sitem D. 
The corresponding field extensions have a chance of being local or graded 
local and are therefore {\it a priori}\/ more interesting than those ``cheap''
ones. 
Among them are the Doplicher-Roberts reconstruction in four dimensions, as
well as some examples in chiral models (related to conformal embeddings) 
which generate only a subset of the superselection sectors of a given model. 
Since we know of no systematic way to construct such extensions, and notably 
we do not have a direct criterium on the canonical triple which ensures 
depth 2, we refrain from discussing this important issue in this report.

The weakness of the Hopf structure reflects the reducibility of the local
inclusions of gauge invariants among the fields. 
According to 3\sitem E, this feature is due to a nontrivial part of the 
quantum symmetry acting innerly on the field algebra. 
The implementing operators lie in the intersection of all local field 
algebras and commute with the observables. 
In the ``regular'' case 5\sitem C with a finite gauge group, they are the 
global implementers of the gauge group, and in the reduced field bundle case 
5\sitem B, they are the source and range projections going along with the 
charged operators \cite{FRS}.

These operators are redundant in order to solve the field problem, in the 
sense that they do not have any effect on the observables. 
Furthermore, they mix up local and global concepts, albeit on the unobservable 
level of charged fields. 
On the other hand, depth 2 and therefore a linear transformation law 
with finite dimensional symmetry tensors require their presence. 
We consider these peculiar field operators as the price to be paid for a 
decent symmetry acting on the charge carrying fields.


\begin{thebibliography}{99}
\itemsep=-.2pc
%
\bibitem{CLZi} K.-H.~Rehren: ``Quantum symmetry associated with braid group
statistics'', in: {\em Quantum Groups}, Proceedings of the ASI Workshop, 
Clausthal 1989, eds.~H.-D.~Doebner {\it et al.}, {\em Lecture Notes in 
Physics} {\bf 370}, pp.\ 318--339, (Springer, 1990).
%
\bibitem{CLZii} K.-H.~Rehren: ``Quantum symmetry associated with braid group
statistics.\ II'', in: {\em Quantum Symmetries}, Proceedings of the ASI 
Workshop, Clausthal 1991, eds.~H.-D.~Doebner {\it et al.}, pp.\ 14--23
(World Scientific, Singapore, 1993).
%
\bibitem{FRS} K.~Fredenhagen, K.-H.~Rehren, B.~Schroer: ``Superselection 
sectors with braid group statistics and exchange algebras.\ I+II'', 
{\em Com\-mun.\ Math.\ Phys.\/} {\bf 125} (1989) 201--226 and 
{\em Rev.\ Math.\ Phys.\/} {\bf Special Issue} (1992) 113--157.
%
\bibitem{BSz} G.~B\"ohm, K.~Szlachanyi: ``A coassociative C*-quantum 
group with non-integral dimensions'', preprint Budapest 1996, q-alg 9509008, 
to appear in {\em Lett.\ Math.\ Phys.}
%
\bibitem{NW} F.~Nill, K.~Szlachanyi, H.-W.~Wiesbrock: ``Weak Hopf algebras 
and reducible Jones inclusions of depth 2'', preprint in preparation, Berlin 
1996;
\newline F.~Nill: unpublished manuscript (1996).
%
\bibitem{Y} T.~Yamanouchi: ``Duality for generalized Kac algebras and a 
characterization of finite group\-oid algebras'', {\em Journ.\ Alg.\/} 
{\bf 163} (1994) 9--50; 
\newline T.~Hayashi: ``Compact quantum groups of face type'',
Publ.\ RIMS {\bf 32} (1996) 351--369.
%
\bibitem{Dr} V.G.~Drinfel'd: ``Quantum groups'', in: Proceedings of the 
Intern.\ Congr.\ of Mathematicians, Berkeley 1986, ed.\ A. Gleason, 
pp.\ 798--820 (Berkeley, 1987).
%
\bibitem{Sch} V.~Schomerus: ``Construction of field algebras with quantum
symmetry from local observables'', {\em Com\-mun.\ Math.\ Phys.\/} {\bf 169} 
(1995) 193--236.
%
\bibitem{O} A.~Ocneanu: ``Quantized groups, string algebras, and Galois 
theory for algebras'', in: {\em Operator Algebras and Applications}, Vol.\ 2, 
eds.~D.E.~Evans {\it et al.}, {\em London Math.\ Soc.\ Lect.\ Notes} 
{\bf 135}, pp.\ 119--172 (Cambridge, 1988).
%
\bibitem{Lo} R.~Longo: ``Index of subfactors and statistics of quantum 
fields.\ I+II'', {\em Com\-mun.\ Math.\ Phys.\/} {\bf 126} (1989) 217--247 
and {\bf 130} (1990) 285--309.
%
\bibitem{LH} R.~Longo: ``A duality for Hopf algebras and for subfactors.\ 
I'', {\em Com\-mun.\ Math.\ Phys.\/} {\bf 159} (1994) 133--150.
%
\bibitem{LR} R.~Longo, K.-H.~Rehren: ``Nets of subfactors'', {\em Rev.\ 
Math.\ Phys.\/} {\bf 7} (1995) 567--597.
%
\bibitem{DHR} S.~Doplicher, R.~Haag, J.E.~Roberts: ``Local observables and 
particle statistics.\ I+II'', {\em Com\-mun.\ Math.\ Phys.\/} {\bf 23} (1971) 
199--230 and {\bf 35} (1974) 49--85.
%
\bibitem{DR} S.~Doplicher, J.E.~Roberts: ``Why there is a field algebra 
with a compact gauge group describing the superselection structure in
particle physics'', {\em Com\-mun.\ Math.\ Phys.\/} {\bf 131} (1990) 51--107.
%
\bibitem{N} F.~Nill, K.-H.~Rehren: unpublished. F.~Nill has already in 1994 
considered the weak Hopf symmetry associated with $\sig\plus$.
\end{thebibliography}
\end{document}